\documentstyle[12pt]{article}
\begin{document}

\title{ Generalized KdV Equation for Fluid Dynamics
and Quantum Algebras}
\author{
A. Ludu$^{1,3}$, R. A. Ionescu$^{2,4}$ and W. Greiner$^1$ \\
\normalsize{$^1$Institut f\"ur Theoretische Physik der 
J W Goethe-Universit\"at,} \\
\normalsize{D-60054, Frankfurt am Main, Germany} \\
\normalsize{$^2$International Centre for 
Theoretical Physics, P.O.Box 586, 34100 Trieste, Italy} \\
}

\date{}
\maketitle
\begin{abstract}
We generalize the non-linear 
one-dimensional equation of a fluid layer for
any depth and length as an infinite order differential equation for the
steady waves. This equation can be written as a q-differential one,
with its general solution  written as a power series expansion with
coefficients satisfying a nonlinear recurrence 
relation. In the limit of 
long and shallow water 
(shallow channels) we reobtain the well known
KdV equation together with its single-soliton solution.
\end{abstract}

\setlength{\baselineskip} {4ex}.

\footnotetext[3]{Permanent Adress: 
Bucharest University, Department of Theoretical Physics, 
P. O. Box MG-5211, Romania;  
E-mail: ludu@th.physik.uni-frankfurt.de}
\footnotetext[4]{Permanent Adress: Institute of Atomic Physics, P.O.Box 
MG-6, Bucharest, Romania; \\ E-mail:remus@ifa.ro}



Short title: {\bf Generalized KdV equation and Quantum Algebras}
\vfill
\eject

\section{Introduction}
\vskip 1truecm

The Korteweg-de Vries 
equation (KdV) was proposed by Korteweg and de Vries
[1] one century ago to 
explain the steady translation water waves observed
in a channel by Scott 
Russel in 1834.  The main property of this equation
consists in the equal occurence of non-linearity and dispersion.  
Some of
its most important solutions, the solitary waves, are single, 
localised
and non-dispersive entities (no dispersion in shape), having also
localised finite energy density. Among these solutions, the solitons 
are
solitary waves with special added requirements concerning their 
behaviour
at infinity ($x \rightarrow \infty $, $t \rightarrow \infty $) and 
having
special properties concerning the scattering between different such
single-solitary solutions. Despite this non-linearity, the KdV 
equation is
an infinite dimensional Hamiltonian system [2-5] and, remarkable, 
when the
KdV solution evolves in time, the eigenvalues of the associated
Sturm-Liouville operator ${\partial }^{2}+\eta $, where $\eta $ is the
solution of the KdV equation, remain constant. Several extended and
complex methods have been developed to study and to solve the KdV
 equation
and other non-linear wave equations (non-linear 
Schr\"odinger
equation, NLS, modified-KdV equation, MKdV, sine-Gordon 
equation, etc.):
inverse scattering method, theoretical group methods, etc. 
(for a review
of some of these techniques one can look into the book by Witham 
[6] or
the paper by Scott et al [7]. 

Due to its properties, the KdV equation was the source of many
applications and results in a large area of non-linear physics (For a
recent review see [8] and the references herein). For example, the KdV
equation and its generalizations (KdV hierarchies, KP equation and its
hierachies and their supersymmetric generalizations [9]) appear in
different fields of physics, from plasma physics [10] to string
 theory
[11], or to describe the dynamics of quantum extended objects [12], 
and
nonlinear phenomena in nuclear physics (non-linear liquid drop model)
[13,14], and the list is far from being exhausted. In all these
applications, the KdV equation (or other non-linear ones) appears 
as the
consequence of certain simplifications of the physical models,
 especially
due to some perturbation techniques involved. 

In the present paper we have used only one of the two 
well-known necessary
conditions for obtaining the KdV equation in shallow water 
channels, i.e.
the smallness of the amplitude of the soliton with respect 
to the depth of
the channel, $h$. This is the single condition which we use, 
and hence, we consider
$h$ to be an arbitrar parameter. In order to obtain the KdV 
equation from
shallow water channel models [3,5], a second condition is imposed,
 i.e. the
depth $h$ of the channel should be smaller than the half-width 
of the
solitary wave. We do not use, in our generalised approach, nor 
this second
condition neither the condition of infinite length of the channel. 
In our
case the length of the channel, $L$, is also considered an arbitrary
parameter. 

We introduce 
this generalisation in connection with the last (above mentioned)
nuclear
physics application. If the height of the 
perturbation is
comparable to the nuclear radius (big clusters formation on the nuclear
surface or symmetric nuclear molecules) it is not appropriate 
to consider
that the nuclear dynamics is that of a shallow fluid layer. 
 Therefore, we should like to study, as a first step of 
our aim, the non-linear dynamics
of a fluid in a channel of arbitrary depth,
in a bounded domain. This different starting
point leads us to a new type of equation which generalizes 
in some sense the KdV equation
and reduces to this in the shallow water case.

As we shall show, this generalized equation (infinite 
order differential
equation) can also be written in terms of a q-differential 
equation. In
this context we introduce here another possible physical
 interpretation
for the q-deformation parameter, in the field of fluid dynamics,
 as the depth of the fluid layer.
We use the connection between the theory of quantum algebras (or 
q-groups, q-deformed algebras), originally initiated by 
Jimbo and Drinfeld [15], and the 
q-differential equations [16].

In the present model we take into account not only the gravity 
but also
the influence of the surface pressure acting on the free surface
 of the
fluid. This implies that the KdV-like structure for the dynamical 
equation
arrises in a lower order (first order in the smallness parameters)
 and one
may not take into account the second order of smallness to 
obtain a KdV
equation. We also note that, despite the general tendency 
followed by 
papers concerning applications of q-groups in physics of first
constructing the deformed algebra and later searching possible
applications, we have obtained a q-differential equation in
 a natural way,
by starting from a traditional physical one-dimensional 
hydrodynamical
model. 

\vskip 2truecm 
\section{The generalized KdV equation}
\vskip1truecm

Let us consider a one-dimensional ideal 
incompressible fluid layer with depth $h$
and constant density $\rho $, in an uniform gravitational field. 
We suppose irrotational motion and consequently the 
field of velocities
is obtained from a potential function $\Phi (x,y,t)$
\begin{eqnarray}
\vec{V}(x,y,t) & = & (u,v)=\nabla \Phi (x,y,t).
\end{eqnarray}
The continuity equation 
for the fluid results in the Laplace equation 
for $\Phi $:
\begin{eqnarray}
\Delta \Phi (x,y,t) & = & 0. 
\end{eqnarray}
The above Laplace equation should be solved
with appropriate
 boundary
conditions for our physical problem. 
We take a
two-dimensional domain: $x \in [x_{0}-L,x_{0}+L]$ (the "horizontal" 
coordinate) and 
$y \in [0, \xi (x,t)]$ (the "vertical" coordinate) where $x_{0}$ 
is an arbitrar parameter, $L$
is an arbitrary lenght and $\xi (x,t)$ is the equation of the free 
surface of the fluid. The boundary conditions on the lateral walls,
$x=x_{0}\pm L$  and on the
bottom $y=0$, consist in  the condition of vanishing of the 
normal velocity
component. The free surface fulfills the kinematic condition
[3]:
\begin{eqnarray}
v\bigg |_{\Sigma } & = & ({\xi }_{t}+{\xi }_{x}u)\bigg |_{\Sigma },
\end{eqnarray}
where  we denote by 
$\Sigma $
the free surface of equation $y=\xi (x,t)$ and the lower index
indicates the derivative.
Eq. (3) expresses the fact that the fluid particles, which belong 
to the 
surface, remain durring time evolution in this surface.
By taking into account the above boundary 
conditions on the lateral walls and on
the bottom, we can write the potential of the velocities 
in the form:
\begin{eqnarray}
\Phi (x,y,t) & = & \sum_{k=0}^{\infty }
\cosh \biggl ({{k\pi y} \over L}\biggr )
\biggl ( 
{\alpha }_{k}(t)
\cos\biggl ({{k\pi x} \over L}\biggr )+
{\beta }_{k}(t)
\sin \biggl ({{k\pi x} \over L}\biggr )
\biggr ),
\end{eqnarray}
where ${\alpha }_{k}$ and ${\beta }_{k}$ are 
time dependent coefficients
fulfiling the condition
\begin{eqnarray}
{\alpha }_{k}
\sin {{k\pi (x_{0}\pm L)} \over {L}} & = &
{\beta }_{k} \cos {{k\pi (x_{0}\pm L)} \over {L}}, \nonumber
\end{eqnarray}
for any positive integer $k$. 
This restriction introduces a special time dependence of 
${\alpha }_{k}$
and ${\beta }_{k}$, i.e. ${{{\alpha }_{k}} \over {{\beta }_{k}}} 
={\gamma }_{k}=const.$, for any $k$, when ${\beta }_{k}\neq 0$.
If ${\beta }_{k}=0$ we simply equate with $0$ the inverse of the above 
fraction. It exists also a dependence
of the ${\alpha }_{k},{\beta }_{k}$ functions of $x_{0}$. This
 special time dependence
does not affect the generality of the potential $\Phi $, but 
only the balance
between the two terms in the RHS of eq.(4). If we fix all 
${\beta }_{k}=0$
we can take arbitrary values for $x_{0}$ and general form for 
${\alpha }_{k}$. In the infinite channel limit,
 $L\rightarrow \infty $,
there are no more restrictions concerning ${\alpha _{k}}$ and 
${\beta }_{k}$
functions, and we can choose  $x_{0}=0$, too. 

We introduce the function:
\begin{eqnarray}
f(x,t) & = & \sum_{k=0}^{\infty }{{k\pi } \over L}
\biggl (
-{\alpha }_{k}(t)
\sin\biggl (k\pi {x \over L}\biggr )
+{\beta }_{k}(t)
\cos \biggl (k\pi {x \over L}
\biggr ) \biggr ), \hfill
\end{eqnarray}
and the velocity field can be written like:
\begin{eqnarray}
u & = & {\Phi }_{x}=\cos(y\partial )f(x,t) \nonumber \\
v & = & {\Phi }_{y}=-\sin(y\partial )f(x,t). \hfill
\end{eqnarray}
where, for simplicity,  the operator $\partial $ represents, 
from now on, the 
partial derivative
with respect to the $x$ coordinate. The equations (5,6) do not 
depend on $L$
and we equally treat the case $L \rightarrow \infty $, when
 we look for
unbounded domains and travelling profiles for the solutions.

In the following we use in the boundary condition (3) the
 velocities
evaluated at $y= \xi (x,t)=h+\eta (x,t)$ in the first 
order in $\eta $,
as we describe the perturbations with small height compared 
to the depth,
but not necessarly with large widths.
Eqs.(6) read:
\begin{eqnarray}
u(x,\xi (x,t),t) & = & \left[ \cos(h\partial )-\eta (x,t)
\partial \sin(h\partial )\right] f(x,t) \nonumber \\
v(x,\xi (x,t),t) & = &  -\left[ \sin(h\partial )+\eta (x,t)
\partial \cos(h\partial )\right] f(x,t). \hfill
\end{eqnarray}

The dynamics of the fluid is described by the Euler
equation at the free surface.
The Euler equation, written on $\Sigma $, results (after writting 
it in terms of the potential, differentiating it with respect to
 $x$ and
imposing the condition $y=\xi (x,t)$) in the form
\begin{eqnarray}
u_{t}+uu_{x}+vv_{x}+g{\eta }_{x}+{1 \over {\rho }}P_{x} 
& = & 0, \hfill
\end{eqnarray}
where $g$ represents the gravitational acceleration and $P$ 
the surface pressure.
Following the same approach used in the calculation of surface 
capilar waves [17],
we have for our one-dimensional case
\begin{eqnarray}
P\bigg |_{\Sigma } & = & {{\sigma } \over {\cal R}}=
{{\sigma {\eta }_{xx}} \over {{(1+{\eta }_{x}^{2})}^{3/2}}}
\simeq -\sigma {\eta }_{xx}, 
\quad \quad \hbox{for small }\eta . \hfill
\end{eqnarray}
where ${\cal R}$
 is the local radius of curvature of the surface (in this 
case,
of the curve $y=\xi (x,t)$) and $\sigma $ is the surface
 pressure coefficient.
In the interior of the fluid the pressure is given by the Euler 
equation.
The nonlinearities 
appear in the dynamics through 
the nonlinear terms in eqs. (3), (7), (8) and (9). 
Consequently, we have a system of two differential equations (3,8)
for the two unknown functions: $f(x,t)$  and $\eta (x,t)$, with 
$u$ and $v$
depending on $\eta $ and $f$ from eqs.(7).
With $f$ and $\eta $ determined and introduced in the
 expressions 
of $u$ 
and $v$, we can finally find the coefficients ${\alpha }_{k}$ and 
${\beta 
}_{k}$.
In the following we treat these equations in the approximation
of small perturbations of the surface $\Sigma $, with respect to 
the depth,
$a=max |{\eta }^{(k)}(x,t)|<<h$, 
where $k=0,...,3$
are  orders of differentiation.

In the linear 
approximation eqs. (3) and (8) become, respectively 
\begin{eqnarray}
-\sin(h\partial )f & = & {\eta }_{t}\nonumber \\
\cos( h\partial )f_{t} & = & -g{\eta }_{x}+{{\sigma } \over 
{\rho }}
{\eta }_{xxx}, \hfill
\end{eqnarray}
and we obtain, by eliminating ${\xi }_{t}$ from eqs.(10),
\begin{eqnarray}
\cos(h\partial ){\eta }_{tt} & = & \sin(h\partial )\biggl (
{{c_{0}^{2}} \over {h^{2}}}
{\eta }_{x}-{{\sigma } \over {\rho }}{\eta }_{xxx} \biggr ), 
\hfill
\end{eqnarray}
which, in the lowest order of approximation in $(h\partial )$
for the sine and cosine functions, 
and in the absence of the surface pressure, 
gives us the familiar wave equation
${\eta }_{tt}=c_{0}^{2}{\eta }_{xx}$, where $c_{0}=\sqrt{gh}$
is the sound wave velocity.
By introducing the solution $\eta =e^{i(kx-\omega t)}$ 
in the linearized eq.(11) 
we obtain a non-linear dispersion relation
$${\omega }^{2}=
c_{0}^{2}k^{2}\biggl (1+{{\sigma } \over {\rho g}}
\biggr ){{\tanh(kh)} \over {kh}}.
$$ 
In the limit of shallow waters we find, for the dispersion 
relations, in 
the case $\sigma =0$ the acoustic waves limit, and, for
$\sigma >> \rho g$ the surface capilar waves limit
${\omega }^{2}={{h\sigma } \over {\rho }}k^{4}
$.
In the absence of the surface pressure ($\sigma =0$) the function
 $f$ is
given, in this linear approximation, at least formally, by:
\begin{eqnarray}
f^{lin} & = & {{c_{0}} \over {h}}\left( 
{{\sin(2h\partial )} \over {2h\partial }}\right)^{-1/2}\eta , 
\hfill
\end{eqnarray}
which in the limit of shallow fluid 
has a particular solution in the form
$$f^{0}(x,t)=
{{c_{0}} \over h}\eta. $$
For the time derivative of $f$ we have, from the second
 equation of eqs.(10),:
\begin{eqnarray}
f^{lin}_{t}(x,t) & = & (\cos (h\partial ))^{-1}
\biggl (-{{c_{0}^{2}} \over 
h}{\eta }_{x}+{{\sigma } \over {\rho }}{\eta }_{xxx}
\biggr ), \hfill
\end{eqnarray}
which in the limit of a shallow fluid  reduces to:
$$f^{0}_{t}(x,t)=-{{c_{0}^{2}} \over h}{\eta }_{x}
+{{\sigma } \over {\rho }}{\eta }_{xxx}$$
Following [3], 
we look for the solution of eqs.(3,8) in the form
\begin{eqnarray}
f & = & {a \over h}c_{0}{\tilde {\eta }}
+\left ({a \over h}
\right )^{2}f_{2}  \hfill \\
f_{t} & = & -c_{0}^{2}{a \over h}(\cos(h\partial ))^{-1}
{\tilde {\eta }}_{x}
+{{a\sigma } \over {\rho }}(\cos(h\partial))^{-1}
{\tilde {\eta }}{xxx}
+{\left({a \over h}\right)}^{2}g_{2}, \hfill
\end{eqnarray}
which represents a sort of perturbation 
technique in $a/h$, where
$\eta=a{\tilde {\eta }}$.
Of course a functional connection exists between the 
perturbations $f_{2}(x,t)$ and $g_{2}(x,t)$.
Eq.(3), in the lowest order in $a/h$, yields:
\begin{eqnarray}
-c_{0}\sin(h\partial ){\tilde {\eta }} & = & 
h {\tilde {\eta }}_{t}+ac_{0}({\tilde {\eta }}_{x}
\cos(h\partial )
{\tilde {\eta }}+{\tilde {\eta }}
\cos (h\partial ){\tilde {\eta }}
_{x}). \hfill
\end{eqnarray}

If we approximate $\sin (h\partial )\simeq h\partial -{1 \over 
6}(h\partial )^{3}$, $\cos (h\partial )\simeq 1-{1 \over 2}
(h\partial )^{2}$, we obtain, 
from eq.(16), the polynomial differential equation:
\begin{eqnarray}
a{\tilde {\eta }}_{t}+2c_{0}{\epsilon }^{2}h{\tilde {\eta }}
{\tilde {\eta }}_{x}+c_{0}\epsilon h{\tilde {\eta }}_{x}
-c_{0}
\epsilon {{h^{3}} \over 6}{\tilde {\eta }}_{xxx}-
{{c_{0}{\epsilon }^{2}h^{3}} \over 2}\biggl (
{\tilde {\eta }}_{x}{\tilde {\eta }}_{xx}+{\tilde {\eta }}
{\tilde 
{\eta }}_{xxx}\biggr ) & = & 0, \hfill
\end{eqnarray}
where $\epsilon ={a \over h}$. The first four terms in eq.(17) 
correspond 
to the zero order approximation terms in $f$, obtained 
from the boundary 
condition at the free surface, in eq.(6.1.15 a) from [3],
i.e. the traditional way of obtaining the KdV equation
in shallow channels.
In this case all the terms are in first and second order in
$\epsilon $.
If we apply the second restriction with respect to the solutions,
i.e. the half-width to be much larger than $h$, we can
 neglect in eq.(17) 
the last paranthesis, and we obtain exactly the 
KdV equation for 
the free surface boundary condition. In other words
we understand the condition
$h\partial $ "small" like $(h\partial )f(x,t)\ll 1$ all over 
the domain of
definition of $f$. This means that the
spatial extension of the
perturbation $f(x,t)$ is large compared to $h$, which is
 exactly the 
case in
which KdV equation arrises from the shallow water model (see 
Chapter 6 in
[3], $h\partial f(x,t)$ of order $\simeq
{h \over{L}}={\delta }\ll 1$).

By using again the approximations given by eq.(7), 
we can write the Euler 
eq.(8) in the form:
\begin{eqnarray}
{\partial }_{t}\Omega f+\Omega f
({\partial }_{t}\Omega f)+
\Omega f(\partial \Omega f) +\omega 
f(\partial \omega f)
& = & -g{\eta }_{x}+{{\sigma }
\over {\rho }}{\eta }_{xxx}, \hfill
\end{eqnarray}
where we used  the notations:
\begin{eqnarray}
\Omega & = & \cos (h\partial )-\eta 
\partial \sin (h\partial ) \nonumber \\
\omega & = & \sin (h\partial )+\eta 
\partial \cos (h\partial ). \hfill
\end{eqnarray}
We note here an interesting property of the operators given
in eq.(19):
\begin{eqnarray}
\Omega +i\omega & = & e^{iy\partial}\bigg | _{\Sigma }+
{\cal O}_{2}(h\partial )+{\cal O}_{2}(h\eta ). \hfill
\end{eqnarray}
Following the same procedure like for 
the free surface boundary condition 
eq.(3), i.e. using the approximation for small $\eta $,
we obtain, from eq.(18):
\begin{eqnarray}
\cos (h \partial )f_{t} & = & -g{\eta }_{x}+
{{\sigma } \over {\rho }}{\eta }_{xxx}, \hfill
\end{eqnarray}
which, in the lowest order in $a/h$ and by using eqs.(14,15),
reduces to an identity.

Before further analysis, we would like to note that in the
shallow water case, following the same addimensional 
notations like in Chapter 6 of [3], i.e. $\epsilon =
{a \over h}$, $\delta ={h \over l}$, where $l$ gives the 
order of magnitude of the half-width of the perturbation 
$\eta $, and introducing the new addimensional parameter
$\alpha =-{{\sigma } \over {gl^{2}\rho }}$, we obtain,
for the Euler equation, the form: 
\begin{eqnarray}
{\tilde {\eta }}_{t^{'}}+
{\tilde {\eta }}_{x^{'}}+{3 \over 2}
{\tilde {\eta }} {\tilde {\eta }}_{x^{'}}+
\alpha {{\epsilon } \over 2}
{\tilde {\eta }}_{x^{'}x^{'}x^{'}} & = & 0, \hfill
\end{eqnarray}
which is again the traditional KdV equation. The primes atached to
the lower indexes represent 
addimensional units, according to [3].
The difference between eq.(22) and the corresponding 
eq.(6.1.15.b) from [3] is given by the inclusion of the
surface pressure effects. If the coefficient $\alpha $
is large enough, one is no more obliged to take into
account second order terms, i.e. ${\delta }^{2}$, to
obtain the KdV equation. Of course, the above approach 
introduces changes in the differential equations which
are involved with higher order perturbations, i.e. $f^{(1)}$
and $f^{(2)}$ from [3] or $f_{2}$, $g_{2}$ in our case.

The reduction of eq.(3) to the KdV equation occures if,
in eq.(16) 
we limit only to  third order of magnitude in $h\partial $:
\begin{eqnarray}
{\eta }_{t}+c_{0}{\eta }_{x}-c_{0}{{h^{2}} \over {6}}{\eta }_{xxx}+
{{2c_{0}} \over h}\eta {\eta }_{x} & = & 0 \nonumber
\end{eqnarray}
In the following we shall investigate this generalised KdV 
equation (GKdV),
obtained from eq.(16), by keeping all terms in $\sin$ and $\cos$, 
namely:
\begin{eqnarray}
{{\eta }_{t}+{{c_{0}}\over h}\sin(h\partial ){\eta } + 
{{c_{0}}\over h}({\eta }_{x}\cos(h\partial )
{\eta }} +\eta \cos (h\partial ){\eta }_{x})
& = & 0 \hfill
\end{eqnarray}
Eqs.(3) and (18), in higher orders in ${a\over h}$ 
and in $h\partial $,
 yield the corresponding differential equations 
for the functions $f_2$ and $g_2$, but we shall study 
here only  eq.(23).

\vskip 2truecm
\section{The q-differential form of the GKdV equation}
\vskip 1truecm

In this Section we limit ourselves to the steady
 translation waves and  
consider only
solutions of the form $\eta (x,t)=\eta (x+Ac_{0}t)=\eta(X)$
 where 
$A\in R$
and $X=x+Ac_0t$. 
Eq.(23) can be written in
 the form:
$$
Ah{\eta }_{X}(X)+{{\eta (X+ih)-\eta (X-ih)} \over {2i}}
+{\eta }_{X}(X){{\eta (X+ih)+\eta (X-ih)} \over 2} 
$$
\begin{eqnarray}
+\eta (X)
{{{\eta }_{X}(X+ih)+{\eta }_{X}(X-ih)
} \over {2}}
& = & 0, \hfill
\end{eqnarray}
if we suppose that $\eta $ is an analytic function in the domain
$Re(z) \in (\infty , \infty )$, $Im(z)\in (-h, h)$ of the 
complex plane. We
study the solutions with a rapid decreas 
at infinity and make a
change of variable: $v=e^{BX}$ for $x \in (-\infty ,0)$ and $
v=e^{-BX}$ for $x \in (0, \infty )$. Here $B$ is a positive real
 constant.
Putting: 
\begin{eqnarray}
\eta (X) & = & -hA+f(v), \nonumber
\end{eqnarray}
we obtain a q-differential equation for the function $f(v)$:
\begin{eqnarray}
f(v){{D_{q}f^{2}_{v}(v)}
\over {D_{q}f_{v}(v)}}
+f_{v}(v){{D_{q}f^{2}(v)} \over {D_{q}f(v)}}+2
{{\sin(Bh)} \over B}D_{q}f(v) & = & 0, \hfill
\end{eqnarray}
where we denoted $q=e^{iBh}$, ($|q|=1$), and the 
operator of
the q-derivative is defined [16]:
\begin{eqnarray}
D_{q}f(v) & = & {{f(qv)-f(q^{-1}v)} \over 
{qv-q^{-1}v}}. \hfill
\end{eqnarray}
We then write the solution of eq.(24) (or eq.(25))
as a power 
series in $v$:
\begin{eqnarray}
f(v) & = & \sum_{n=0}^{\infty }a_{n}v^{n}, \hfill
\end{eqnarray}
and we must put $a_{0}=hA$ to have $\lim_{x \rightarrow 
\pm \infty}\eta
(x)=0$. Equation (25) results in a recurrence non-linear
 relation for the
coefficients $a_{n}$:
$$
\left[ Ahk+{1 \over B}\sin(Bhk)\right]a_{k} 
$$
\begin{eqnarray}
& = & -
\sum_{n=1}^{k-1}n\biggl (\cos\left( Bh(k-n)\right )
+\cos (Bh(k-1))\biggr )a_{n}a_{k-n}. 
\hfill
\end{eqnarray}
By taking $k=1$ in the above relation we obtain $a_{1}
\left [Ah+
{1 \over B}
\sin(Bh) \right ]=0$. 
Without any loss of generality, due to the
arbitraryness of $B$ we can write:
\begin{eqnarray}
A & = & -{{\sin(Bh)} \over {Bh}}. \hfill
\end{eqnarray}
This relation fixes the velocity of the envelope of the 
perturbation if its
asymptotic behaviour is fulfilled. To have $A \neq 0$ we 
need $Bh\neq k\pi $
for $k$ integer. In this condition $a_{1}$ is arbitrary 
at present and, by
writing $a_{k}={\alpha }_{k}a_{1}^{k}$ we have
 ${\alpha }_{1}=1$ and the
recurrence relation:
\begin{eqnarray}
{\alpha }_{k} & = & {{2B\cos {{Bh(n-1)} \over 2}
} \over 
{k\sin(Bh)-\sin(kBh)}}\sum_{n=1}^{k-1}
n
\cos {{Bh(2k-n-1)} \over 2}
{\alpha }_{n}{\alpha }_{k-n}, 
\hfill
\end{eqnarray}
for $k \ge 2$. This recurrence relation gives the $k$ 
coefficient if we know
the first $k-1$ ones. 

To obtain a smooth behaviour of the solution $\eta
(X)$
in $X=0$, i.e. a continuity of its derivative, we must introduce 
the
condition:
\begin{eqnarray}
f_{v}(1) = \sum_{n=1}^{\infty }n{\alpha }_{n}
a_{1}^{n-1} & = & 0, 
\hfill
\end{eqnarray}
or the derivative of the power series $f(v)$ with the coefficients
 given by
eq.(28) is zero in $z\in R, z=a_{1}$. This condition fixes $a_{1}$.

We study now a limiting case of the relation (28), by replacing 
$\sin$ 
and $\cos$
expressions with their lowest nonvanishing terms in their power 
expansions.
Thus, we have:
\begin{eqnarray}
{\alpha }_{k} & = & {6 \over {
B^{2}h^{3}k(k^{2}-1)
}}
\sum_{n=1}^{k-1}n{\alpha }_{n}{\alpha }_{k-n}. \hfill
\end{eqnarray}
It is a straightforward exercise to prove that: 
\begin{eqnarray}
  {\alpha }_{k} & = &
\biggl ( {1 \over {2B^{2}h^{3}}}\biggr )^{k-1}k  \hfill
\end{eqnarray}
is the solution of the above recurrence relation. This can be
 done by using
the mathematical induction and by taking into account the relations:
\begin{eqnarray}
\sum_{n=1}^{k-1}n^{2} & = & {{k(k-1)(2k-1)} \over 6} 
\nonumber \\
\sum_{n=1}^{k-1}n^{3} & = & \left({{k(k-1)}\over 2}\right)^2. 
\nonumber
\end{eqnarray}
We can write the power expansion:
\begin{eqnarray}
g(z) & = & \sum_{k=1}^{\infty }k\biggl ( 
{1 \over {2B^{2}h^{3}}} 
\biggr )^{k-1} z^{k}. \hfill
\end{eqnarray}
having the radius of convergence ${\cal R}=2B^{2}h^{3}$ (due to 
the 
Cauchy-Hadamard criterium).
The function $g(z)$ can be written in the form
\begin{eqnarray}
g(z) & = & z \biggl ({1 \over {1-{z \over {2B^{2}h^{3}}}}} 
\biggr )_{z}
2B^{2}h^{3}=-{z \over {
\left( 1-{z \over {2B^{2}h^{3}}} \right)^{2}
}} \hfill
\end{eqnarray}
The condition (31) results in $a_{1}=-2B^{2}h^{3}$ and 
\begin{eqnarray}
{\alpha }_{k} & = & k\biggl ( {1 \over {2B^{2}h^{3}}}
 \biggr )^{k-1}
(-2B^{2}h^{3})^{k}=2B^{2}h^{3}(-1)^{k}k, \hfill
\end{eqnarray}
which gives:
\begin{eqnarray}
\eta (X) & = & 2B^{2}h^{3}\sum_{k=1}^{\infty }k
\biggl ( -e^{-B|X|}
\biggr )^{k}= \nonumber \\
& & 2B^{2}h^{3}
{{e^{-B|X|}} \over {\left( 1+e^{-B|X|}\right)^{2}}} 
= \nonumber \\
& & {{B^{2}h^{3}} \over 2}{1 \over {\left(
\cosh(BX/2)\right)^{2}}} .
 \hfill
\end{eqnarray}
As one expects, the above solution is exactly the single-soliton 
solution
of the
KdV equation and it was indeed obtained by assuming $h$ small 
in the
recurrence relation (28).
We stress at this point that the GKdV equation has two general 
features,
expressed in the reduction of eq.(16) to KdV equation and in 
eq.(23):
in the limit $h\partial $ "small", both the differential 
equation GKdV
and one of its solutions, $\eta (X)$, go versus the KdV equation 
and its
single - soliton solution, correspondingly.

In the general case we do not obtain a simple solution like the 
relation 
(33) but we have all the necessary informations from the 
relations (28).
It seems that the power series $g(z)$ with the coefficients given
 by
the recurrence relations (28) has a nonvanishing radius of
 convergence. The 
problem of the
existence of a real point $z_{0}$ in the disk of convergence or on its 
frontier in which $g^{'}(z_{0})=0$ needs 
further studies. For the KdV equation this point is on the 
frontier of the
disk of convergence of the power series with coefficients
 given by the 
relations (28).
We mention that $g^{'}(z_{0})=0$ means that the function 
is not univalent (
it is not injective in the vicinity of this point).

An interesting comment concernig an unexpected 
connection is appropriate: Bieberbach
formulated in 1916 the following conjecture: 
if $f(z)= z+a_2 z^2 +...$ is
analytic and univalent in the unit disk, then $|a_n|\le n$ 
for all $n$, 
with equality occuring only for rotations of the Koebe 
function 
$k(z)={{z}\over{(1-z)^2}}$ [18]. The form of our solution 
for the 
KdV equation, eq. (34), is exactly the Koebe function up 
to a scaling of 
the $z$ variable. This conjecture was demonstrated in 1986
 by de Branges
[19]. This fact suggests that one could obtain information 
about 
the possible zeros of the $g'(z)$ if it is possible to
 obtain estimates 
for the convergence radius of the power series with the coefficients 
given by (30).

The change of the variable which we used to obtain the non-linear 
q-differential equation (25) can be used in the time-dependent
equation obtained from
eq.(23). If $\eta(x,t)=f(v,t)$ where $v=e^{Bx}$ or
 $v=e^{-Bx}$, then
we obtain:
\begin{eqnarray}
2{{h} \over {c_{0}}}{1\over v}f_{t}^{\mp}(v,t)\pm B
f_{v}^{\mp}(v){{
D_{q}f^{\mp 2}(v,t)} \over 
{D_{q}f^{\mp}(v,t)}} \nonumber \\
\pm Bf^{\mp }(v){{D_{q}f_{v}^{\mp 2}(v,t)
} \over {D_{q}f_{v}^{\mp }(v,t)}}
\pm 2{\sin(Bh)} 
D_{q}f^{\mp}(v,t) & = & 0, 
\end{eqnarray}
where $\pm$ reffers to the positive or negative real semiaxis.
If we try a solution of the form $f^{\pm}(v,t)=\sum_{k=0}^{\infty 
}a_{k}^{\pm}v^{k}$,
then we obtain the following infinite system of differential 
equations:
$$
{h \over {c_{0}}}{{d}\over{dt}}a_{k}^{\mp}(t)  
={\pm} \sin(kBh) 
a_{k}^{\mp}(t)
$$
\begin{eqnarray}
& \mp & 2B\cos{{Bh(n-1)} \over 2} 
\sum_{n=0}^{k}na_{n}^{\mp}(t)a_{k-n}^{\mp}
(t)\cos \biggl (
{{(2k-n-1)Bh} \over 2}
\biggr ). 
\end{eqnarray}
We note that if $a_0(t_0) = 0 $ then this relation 
is preserved 
in 
time. The above system can be resolved step by step as
 the equation 
for $a_n$ involves only the coefficients $a_k$ with $k\le n$. 
However, 
the smooth behaviour of the solution in $x=0$ 
at every moment is a 
nontrivial problem. 
\vskip 1.5truecm

\section{Conclusions}
\vskip 1.5truecm

We have shown that the KdV equation describing the shallow 
liquids can be
generalised for any depths and lengths.
Consequently
the influence of the layer depth could be included in the present 
studies on
different fields, like clusterization on the nuclear surface in nuclear
physics.

The present paper shows two important results: a generalization 
of the KdV
equation starting from a physical model and an embeding of
 this non-linear
equation into a q-deformed formalism.

By using a non-linear hydrodynamic approach for a fluid layer 
we obtained 
a dynamical differential equation, of infinite order, wich 
generalizes the
KdV equation for shallow water. Contrary to the later one, 
the generalised
equation is valid for any depth and length of the fluid layer. 
We succeded to write this equation in a q-differential 
formalism, and
consequently we have obtained a non-linear recurrence relation 
for the 
coefficients of its general solution. 
We stress the importance of the introduction of the surface
pressure term, eqs.(8,9), which provides the necessary 
dispersion term, i.e. ${\eta }_{xxx}$.
This term is esential in two respects:
first, it introduces the dispersion in a smaller order than
in the traditional case [3], and second,
it represents the single term responsible for dispersion
in the case of cylindrical geometry [13,14]. 
Both the generalized
 KdV equation 
and its formal solution reach the KdV system in the
 shallow water 
limiting case, which gives confidence in the present 
approach.
We conjecture that exists a deep connection between 
non-linear
 differential
equations, infinite order linear differential equations 
(which are connected 
with
finite diference equations) and q-differential equations 
and their q-deformed
symmetries. In this sens, in the present paper we succeded 
to show that by
starting from a one dimensional model for an ideal 
incompressible
irrotational fluid layer, a more general differential equation, 
than the KdV
equation, was obtained.
The equation which we proposed  not only 
generalizes the KdV equation, 
but also can be written as a q-deformed nonlinear 
equation with a 
deformation parameter of modulus one. In our case the
 deformation parameter 
$q$ is put into
relation with the depth of the fluid. 
To our knowledge this one of the few 
cases when one can add a physical
interpretation to this parameter.

We stress that the result of the  Section 3, i.e. the realisation
 of the
GKdV equation as a q-deformed one, could be the starting point  
to search for more interesting symmetries. Perhaps, there 
is a more general 
connection with the q-deformed algebras. 
It would be interesting to interpret the
q-GKdV eqs.(24,25) as the Casimir element of a certain 
q-deformed algebra.
This posibility could open the way
to construct invariants of this equation and to compare them, in the
limiting case of shallow waters, with the KdV invariants. 
Such work is in progress.

We expect many possible application of such a formalism, 
especially in the
field of nuclear clusterization. If the 
GKdV equation arrises from a
Lagrangean formalism or from a Hamiltonian formalism, then, 
one can apply this
result in physical models which involve exotic shapes (cluster decay,
spontaneous fission, fragmentation of nuclei, atomic clusters, etc).

\vskip 1truecm

{\bf Acknowledgements:} One of the authors (R.A.I.) would like 
to thank Professor Abdus Salam, the IAEA and UNESCO for 
hospitality at 
the International Centre for Theoretical Physics, Trieste.

\end{document}